\documentclass[10pt,notitlepage,twoside,a4paper]{amsart}

\usepackage{amsmath,amssymb,enumerate}

\usepackage{epsfig,fancyhdr,color}%,showkeys,amsmidx

\usepackage{amssymb}
\usepackage{amsmath,amsthm}    
\usepackage{latexsym}
\usepackage{amscd}

\usepackage{latexsym}
\usepackage{amscd}

\usepackage[all]{xy} 

\usepackage[mathcal]{eucal}
\usepackage{pgf,tikz}

\definecolor{NoteColor}{rgb}{1,0,0}

\newtheorem{theorem}{\rm\bf Theorem}
\newtheorem{proposition}{\rm\bf Proposition}[section]

\theoremstyle{definition}
\newtheorem{definition}[proposition]{\rm\bf Definition}

\theoremstyle{remark}

\def\interieur#1{\mathord{\mathop{\kern 0pt #1}\limits^\circ}}

\title{On three early papers by Herbert Busemann}

\author{Athanase Papadopoulos} 
\address{Institut de Recherche Math\'ematique Avanc\'ee,
Universit{\'e} de Strasbourg and CNRS, 
7 rue Ren\'e Descartes, 
 67084 Strasbourg Cedex, France.}
\email{papadop@math.unistra.fr}

\author{Marc Troyanov} 
\address{Section de Math{\'e}matiques,  
\'Ecole Polytechnique F{\'e}derale de Lausanne, station 8,
1015 Lausanne - Switzerland} 
\email{marc.troyanov@epfl.ch}

\begin{document}

\begin{abstract}
This paper is a commentary and a reading guide to three papers by Herbert Busemann,   
\emph{\"Uber die Geometrien, in denen die  ``Kreise mit unendlichem Radius''   die k\"urzesten Linien sind.}   (On the geometries where circles of infinite radius are the shortest lines) (1932), \emph{Paschsches Axiom und Zweidimensionalit\"at}, (Pasch's Axiom and Two--Dimensionality) (1933) and \emph{\"Uber  R\"aume mit konvexen Kugeln und Parallelenaxiom. }
 (On spaces with convex spheres and the parallel postulate) (1933). These are the first papers that Busemann wrote on the foundations of geometry and the axiomatic characterization of Minkowski spaces (finite-dimensional normed spaces). The subject of these papers followed Busemann for the rest of his life, and the three papers already contain several ideas and techniques that he developed later on, in his work on the subject which lasted several decades. The three papers were translated into English by Annette A'Campo. These translations, together with the final version of present commentary, will be part of the forthcoming edition of Busemann's \emph{Collected papers} edition.

\end{abstract}
\maketitle

%_______________


\medskip


\medskip
% \tableofcontents   

\section{Introduction}  \label{intro}
  
The three papers \cite{Busemann1932b,Busemann1933a, Busemann1933}  by Busemann with which we are concerned in this note are among the first papers that he wrote.\footnote{Busemann published a paper on analysis in 1930, \emph{Die Vollst\"andigkeit der Minimalfolgen von Eigenwertproblemen}, Nachrichten G\"ottingen 1930, 295-307. 
It seems he did not publish any paper in 1931.} The first one was published in 1932 and the two others in 1933. These papers concern the foundations of geometry, a field whose main object is to investigate the various existing geometries, in general from the axiomatic point of view. This concerns in particular the three ``classical''  geometries: Euclidean, hyperbolic and spherical, but projective and affine geometries were also extensively studied, since the 19th century. Furthermore, the axiomatization process, when Busemann started his work, had already been extended to other fields than geometry. One may mention arithmetic, but we also recall that  Hilbert's Problem VI, from his famous 1900 twenty-three problems, asks for the axiomatization of mathematical physics.
 
 The axiomatic method has a long history which starts in  Greek Antiquity\footnote{It is well known that the Greeks started the axiomatic foundation of Euclidean geometry, but they also worked out the foundations of spherical geometry; this was done mainly by Theodosius and Menelaus. They also approached the foundations of  projective geometry; one may think here of the work of Apollonius.} with a culmination in Hilbert's \emph{Grundlagen der Geometrie} \cite{Hilbert-Grundlagen}. The first edition  of this treatise appeared in 1899, and the seventh revised edition, to which Busemann refers, in 1930. 

It is also worth mentioning that during the three decades that preceded Busemann's work, set theory and topology, which were closely related fields, experienced a rapid development, and several mathematicians, including Busemann himself, were familiar with them.\footnote{We point out that Busemann, during his studies, followed the topology courses of P. S. Alexandroff, who visited G\"ottingen regularly. See \cite{Goenner}.}  Busemann had also at his disposal the theory of metric spaces.  We recall in this respect that the abstract notion of metric space as we understand it today first appeared in 1906, in Fr\'echet's thesis \cite{Frechet}. At the time where Busemann started his work, this theory had already witnessed important developments, in particular by K. Menger who worked in Vienna. Busemann refers at several places to the work of Menger.

The papers  \cite{Busemann1932b,Busemann1933a, Busemann1933}   by Busemann are important because the questions that Busemann considers and the ideas he expresses there are 
 among those that followed him during the rest of his life.
At the time where Busemann wrote them, he was working as an assistant (without salary) at the University of G\"ottingen, having obtained his doctorate in 1931 under the supervision of R.  Courant. 
Hilbert, to whom we shall often refer in this note, was still teaching at that university, mostly on philosophy of mathematics, although he had officially retired in 1930.

The specific problem in which Busemann was interested in these three papers is that of characterizing  Minkowski geometries. We recall that a Minkowski space is defined by a metric  $\rho$ on a finite-dimensional real vector space such that for any quadruple $X,Y,X',Y'$, if  $(Y'-X') = t (Y-X)$, then $\rho(X',Y') = |t| \rho(X,Y)$,
see  \cite[chap 1, \S 1]{Minkowski-Zahlen}  .
 These spaces were   introduced by Hermann Minkowski in the book \cite{Minkowski-Zahlen}, to which Busemann refers. These spaces play a major role in Busemann's subsequent work. They are the infinitesimal Finsler geometries, that is, they play, in Finsler geometry, the role played by Euclidean spaces in Riemannian geometry.  Developing a theory of Finsler spaces from a purely metric point of view, that is, without differentiability, was one of  Busemann's favorite topics. Another reason for which Minkowski spaces are important in Busemann's work is that they constitute one of the two main classes of examples of metric spaces satisfying the conditions of Hilbert's Problem IV.\footnote{The other class is that of Hilbert metrics.} This problem, from the list of problems that Hilbert addressed at the Paris 1900 International Congress of Mathematicians, remained for several decades one of Busemann's favorite problems. In Hilbert's words, this is the \emph{problem of the straight line as the shortest distance between two points}. This is usually interpreted as the question of characterizing the metrics on subsets of Euclidean (or projective) space such that the Euclidean (respectively projective) straight lines are the shortest lines. As a matter of fact, in his comments on that problem, Hilbert says much more. In particular he mentions twice Minkowski spaces.   First, he considers them as one of the two  main examples of spaces satisfying his problem, and then he asks for an axiomatic characterization of these spaces. Let us quote his words, from the English translation of his fourth problem \cite{Hilbert-Problems}:
 \begin{quote}\small
 One finds that such a geometry really exists and is no other than that which Minkowski constructed in his book, Geometrie der Zahlen,\footnote{[Hilbert's footnote] Leipzig, 1896.} and made the basis of his arithmetical investigations. Minkowski's is therefore also a geometry standing next to the ordinary Euclidean geometry; it is essentially characterized by the following stipulations:
 \begin{enumerate}[1.]
  \item The points which are at equal distances from a fixed point $O$ lie on a convex closed surface of the ordinary Euclidean space with $O$ as a center.
  \item Two segments\footnote{In this paper, the word ``segment" is used several times, in more than one sense. In each
case the meaning should be clear from the context.} are said to be equal when one can be carried to the other by a translation of the ordinary Euclidean space.
\end{enumerate}
In Minkowski's geometry the axiom of parallels also holds under the conditions that balls are strictly convex. 
 \end{quote} 
 
 We note that in his formulation of the second property, Hilbert makes use of translations of the space. In Busemann's axiomatic characterization of Minkowski spaces, which is the subject of the two papers \cite{Busemann1932b} and  \cite{Busemann1933}, there is no mention of motion. 
 We remind the reader that the question of including or not the motions in the axioms of geometry has been a delicate issue, for several centuries. It originates in Aristotle, who considered that motion should be avoided in the axioms of geometry, because it pertains to physics rather than mathematics, cf. e.g. 
Aristotle's Metaphysics  A8, 989, b32-33  \cite{Aristotle-Metaphysics}.
Euclid, in his list of axioms, followed Aristotle's point of view.  The question of including or not motions in the foundations of geometry has been thoroughly discussed by several commentators of Euclid, in particular the Arabs; see e.g. \cite{RH}, \cite{R} and the references there.

The question of describing axioms for Minkowski geometries remained central for Busemann. It is investigated in his  later papers, e.g. \cite{Busemann-Minkowski1950,BP1979}, and in his books \cite{Busemann1955} and \cite{Busemann1970}.  Besides the question of characterizing Minkowski spaces, Busemann made an extensive investigation of the geometry of these spaces, cf. \cite{B-Minkowski1969,BP1979}.
There is more about these results in the paper \cite{Berestovskii} in the present volume.

 We shall comment more on Hilbert's Problem IV. In particular, we shall  point out that several directions that Hilbert indicated, in the comments that accompany this problem, concerning the solution.  Busemann remained interested in that problem during several decades; cf. \cite{2012-Hilbert1} for a review of his work on this problem.
He started by following Hilbert's line of thought, but he then brought a certain number of fundamental ideas that were completely novel. His program combines the then recent developments of metric geometry, set theory and topology, together with a knowledge of the classical
bases of axiomatic and projective geometry. At several occasions, Busemann commented on his metric approach to the axioms. In the first two paragraphs of  \cite[Chap. 3]{Busemann1942}, he writes:

\smallskip

\begin{quote}\small 
 The axiom that there is only one straight line through two given points is fundamental in the foundation of geometry, in projective geometry 
 and in non-euclidean geometry.  At the stage  where this axiom is formulated a metric is generally not yet defined; therefore it is meaningless to speak of 
 the straight line as a geodesic. 
 
  The early geometers frequently tried to introduce the straight lines as  shortest connections,\footnote{For instance, Legendre in his Elements de g\'eom\'etrie
 gives the following definition: \emph{La ligne droite est une ligne ind\'efinie qui est le plus court chemin entre deux quelconques de ses points  (The straight line
 is an undefined line that is the shortest path between any two of its points)}.} but without success since there was no abstract definition  of a metric.
\end{quote}
  
Hilbert suggested in his comments on Problem IV the inclusion of metric elements as part of the axiomatics. More precisely, he asks whether the triangle inequality could be 
introduced as an axiom in the foundation of geometry. He writes in \cite{Hilbert-Problems}: 

\smallskip

\begin{quote}\small 
We are asking, then, for a geometry in which all the axioms of ordinary Euclidean geometry hold, and in particular all the congruence axioms except the one of the congruence of triangles (or all except the theorem of the equality of the base angles in the isosceles triangle), and, in which, besides, the proposition that \emph{in every triangle the sum of two sides is greater than the third} is assumed as a particular axiom. 
\end{quote} \label{quote.hilbert}

Regarding this passage, let us recall that the triangle inequality is a \emph{theorem} in Euclidean geometry; see the end of the Appendix for a short discussion.

 We mentioned that the paper \cite{Busemann1933a}  concerns the  axioms of two-dimensional geometry. We note in this respect that it may happen that in some  results of Busemann, the case of dimension two requires a special statement (cf. for example the results of the paper \cite{Busemann1933} that are recalled in \S \ref{s:1933} below). Since we are talking about dimension, we note right away that there exist axiomatic definitions of dimension (in the style of Hilbert's \emph{Foundations}), but that Busemann in his work uses the notion of topological dimension, for which he refers to Alexandroff \cite{Alexandroff}.

The plan of the rest of this note is the following.  

In \S \ref{s:notions}, we introduce some basic notions that are used by Busemann, in particular the notion of line space.
 
 In \S \ref{s:Pasch}, we comment on the content of  the paper  \cite{Busemann1933a} where Busemann proves
 that a line space is 2-dimensional if and only if Pasch's Axiom holds.

 In \S \ref{s:1932b}, we comment on the content of Busemann's paper \cite{Busemann1932b}. The main result in this paper is a characterization of Minkowski spaces by the fact that horospheres, with the induced metric, are themselves line spaces.

 In \S \ref{s:1933}, we comment on the content of Busemann's paper \cite{Busemann1933}. In this paper, a characterization of Minkowski spaces is given in terms of a parallelism relation satisfied by the space together with a convexity property of its spheres. This property of convexity of spheres will play a prominent role in the later works of Busemann. 
 
 In an appendix at the end of the present note, we have included, for the convenience of the reader, the list of Hilbert's axioms of Euclidean geometry to which Busemann refers.

%____
\section{The notion of line space}   \label{s:notions}

In the three papers, Busemann uses the following notion that he attributes to Menger.  

\begin{definition}    A \emph{line space} is a set equipped with a non-negative distance function $\rho(X,Y)$, called a metric, satisfying the following 6 axioms.
 
 \smallskip

\begin{enumerate}
\item[\rm{i}.] $\rho(X,Y)=\rho(Y,X)$.
\item[\rm{ii}.] $\rho(X,Y)=0$ if and only if $X=Y$.
\item[\rm{iii}.]   $\rho(X,Y)+\rho(Y,X)\geq \rho(X,Z)$.
\item[\rm{iv}.] Every sequence which is bounded with respect to this metric has an accumulation point.
\item[\rm{v}.]  For every pair of points $(X,Y)$ there exists a unique point $U$ such that 
  \[\rho(X,U)=\rho(U,Y) \quad  and  \quad  \rho(X,U)+\rho(U,Y)=\rho(X,Y).   \]
\item[\rm{vi}.]  For every pair of points $(X,Y)$ there are exactly two points $V$ and $W$ such that \\
 \[\rho(X,Y)=\rho(Y,V)=\rho(W,X),\]
 \[\rho(X,Y)+\rho(Y,V)=\rho(X,V),\]
 \[\rho(Y,X)+\rho(X,W)=\rho(Y,W).\]
\end{enumerate}
 
\begin{tikzpicture}[line cap=round,line join=round,x=1.0cm,y=1.0cm]
\clip(-4.5,-0.4) rectangle (6,1);
\draw [domain=-2.56:5.18] plot(\x,{(-0-0*\x)/3});
\fill (0,0) circle (1.5pt);
\draw(0.1,0.25) node {$X$};
\fill (2,0) circle (1.5pt);
\draw(2.09,0.25) node {$Y$};
\fill  (1,0) circle (1.5pt);
\draw (1.09,0.25) node {$U$};
\fill  (4,0) circle (1.5pt);
\draw (4.1,0.25) node {$V$};
\fill (-2,0) circle (1.5pt);
\draw(-1.91,0.25) node {$W$};
\end{tikzpicture}
\end{definition}
 
 The first three axioms are the usual axioms of a metric space. Axiom iv is a formulation of the Bolzano-Weiertrass property for sequences in $\mathbb{R}^n$. It seems that the introduction of this statement in the axioms of a geometry is a novelty.
 
 Axiom v is about the existence and uniqueness of \emph{midpoints} and Axiom vi is about the existence and uniqueness of \emph{doubled points}.

Busemann notes that, using  techniques from metric geometry as in Menger's paper \cite{Menger1928}, one may prove that 
 \emph{any pair of distinct points $X$ and $Y$ in a line space belongs to a unique geodesic  that is isometric to $\mathbb{R}$.} He attributes this fact to Biedermann. 
Busemann also notes that in a line space, Hilbert's linear axioms are satisfied:\footnote{The numbering here refers to that of Hilbert's axioms in the appendix of this article.}  I, 1,2,3 (first half); II, 1,2,3, III, 1,2,3 and  V.1.

In his later works, Busemann calls such a space  a \emph{straight line space} (e.g. in \cite{Busemann1942}) or simply  a \emph{straight space} 
(e.g. in \cite{Busemann1955}).

\section{The paper on Pasch's axiom and 2-dimensionality}
\label{s:Pasch}
In the short paper \cite{Busemann1933a}, Busemann recalls  Pasch's axiom\footnote{It may be useful to recall that Pasch's axiom was known and used by Euclid's commentators long before Pasch (1843--1930); see e.g.  \cite{RH}.}  which states that for any three non-collinear points $A,B,C$ in a plane, if a line meets the segment $AB$, then it meets the union of the segments $AC$ and $CB$. He proves that under some topological conditions, Pasch's axiom is equivalent to the fact that the space is homeomorphic to the plane. More precisely, Busemann shows that if the lines (in the sense of axiomatic geometry) are the geodesics for a metric satisfying the axioms of a line space as defined in \S\,\ref{s:notions}, then \emph{the fact that Pasch's axiom holds is equivalent to the fact that the space is homeomorphic to the plane}. The proof is sketched in \cite{Troyanov}.

Note that in a line space, one can easily define the congruence of segments by stating that $AB \equiv A'B'$ if $\rho(A,B) = \rho(A',B')$,
and this congruence satisfies Hilbert's axioms (III.1) to (III.3).  

In the last part of the paper, Busemann proposes the following notion of 
angle congruence in a 2-dimensional line space. Let $g$ and $h$ be two rays 
with common origin $O$ and $g',h'$ two other rays with common origin $O'$. Consider the points $G$ and $H$ on $g$ and $h$ respectively satisfying $\rho(O,G)=\rho(O,H) = 1$ and similarly the point  $G'$ and $H'$ on $g'$,$h'$ such that $\rho(O',G')=\rho(O',H') = 1$. Then the angles $(g,h)$ and $(g',h')$ are said to be congruent if $\rho(G,H) = \rho(G',H')$.
This angle congruence is  clearly an equivalence relation and Hilbert's axiom (III.4) is satisfied.

In general however, Axiom (III.5)  is not satisfied and therefore the first theorem on triangle congruence in Hilbert's \emph{Foundations}
fails for a general line space. We recall that this theorem states the following (sometimes called the  SAS -- for side-angle-side -- 
triangle congruence criterion).

\smallskip 

\emph{If $ABC$ and $A'B'C'$ are two triangles such that $A'B' \equiv AB$,  $A'C' \equiv AC$ and $\widehat{B'A'C'} \equiv \widehat{BAC}$,
then triangles $ABC$ and $A'B'C'$ are congruent, in particular we also have $B'C' \equiv BC$.}

\smallskip 

Busemann concludes his paper by stating that \emph{it is well known that one cannot say more in general, even with a quadratic metric}.
Busemann probably means here that in a two-dimensional Riemannian manifold (say complete, 
simply connected and without conjugate points), although the notion of angle is well founded, the theorem on triangles congruence 
generally fails.

%_______________________

\section{Geometries where Circles with Infinite Radius are the Shortest Curves} \label{s:1932b}

Let us now discuss the content of the  paper \emph{On the Geometries
where Circles with Infinite Radius are the Shortest Curves} \cite{Busemann1932b}.  
Busemann starts with the following definition:

\begin{definition}  
 A  line space   satisfies the \emph{limit circle axiom} if in addition to (i)--(vi)
 above, it also satisfies the following condition:
 \begin{enumerate}
\item[\rm{vii}.]  For any pair of distinct
  points $X,Y$ and any sequence $\{ P_n \} $ such that $\rho(
  P_n, X) \rightarrow \infty$ and $( \rho( P_n, X) - \rho( P_n, Y))
  \rightarrow 0$ we also have
  \begin{eqnarray*}  %\label{cond.LimitCircle}
    ( \rho( P_n, X) - \rho( P_n, U)) \rightarrow 0 & \text{and} & ( \rho(
    P_n, X) - \rho( P_n, V)) \rightarrow 0,
\end{eqnarray*}
  where $U$ is the midpoint of $X$ and $Y$ and $V$ is the doubling point, that is
  $Y$ is the midpoint of $X$ and $V$. 
\end{enumerate}

\end{definition}

\textbf{Remaks}. 1.) Condition  vii  is equivalent to the requirement
that
 $$( \rho( P_n, X) - \rho( P_n, Z)) \rightarrow 0$$ 
 for any point $Z$ on
the line through $X$ and $Y$. In Euclidean geometry, this holds if
the sequence $P_n $ diverges to infinity on a line othogonal to the line
through $X$ and $Y$.

2.) This condition also  implies   that the horospheres in  our line space  are
totally geodesic.

The goal of the paper is to prove the following

\begin{theorem}
  A 3-dimensional line space  satisfying the limit circle axiom is isometric to a
  Minkowski space. Furthermore, the metric balls  are strictly convex and the spheres have a unique
  tangent plane at any point.
\end{theorem}

Busemann also proves the converse: every  Minkowski
space with the above properties is a line space satisfying the limit circle axiom.

Let us describe in some detail the content of the paper. The article contains 5 sections comprising a total of 41 
facts which are either technical lemmas or propositions having some interest in themselves.

In Section 1, Busemann describes the segments and  lines in a line space and spells their basic properties.
The results can be summarized by saying that the usual incidence axioms for lines such as those listed in Hilbert's 
\emph{Grundlage der Geometrie} are satisfied.

\smallskip

Let us note Fact 3,  which states that the join of a point $P$ and a segment $AB$ is homeomorphic to a triangle 
in $\mathbb{R}^2$ if $A,B$ and $P$ are not on a line.

In Section 2, Busemann introduces the notion of \emph{limit sphere} in a general line space. 
By definition, a limit sphere is the limit of a sequence of spheres through a point $P$ with centers $Q_n$
lying on a fixed ray and such that $\rho(P,Q_n) \to \infty$.
The limit spheres are also called \emph{horoshperes}.\footnote{It is also worth recalling that the notion of limit sphere was introduced by Lobachevsky, in the setting of his work 
on hyperbolic geometry. We refer the reader to the English edition of the Pangeometry \cite{Pangeometry} (p. 8). We also note that Lobachevsky, in some other memoirs, used the word \emph{horisphere} (see e.g. \cite{Loba-Geometrische}), which was later on transformed into \emph{horosphere}. Lobachevsky proved that in hyperbolic space, horospheres are equipped with a Euclidean structure. His proof is based on the fact that Euclid's axioms are satisfied on horospheres, for a appropriate notion of line, angle, etc.}

Note that Busemann introduced the notion of horosphere in a general metric space, as a level set of a \emph{Busemann functions}
in  \cite[p. 102 and p.132]{Busemann1955}, see also \cite{Papadopoulos2005}. 
Recall that for any geodesic
ray $\gamma$, the associated Busemann function $\beta_{\gamma}$
is defined as
$$
 \beta_{\gamma} (X) = \lim_{t \to \infty} (\rho(X,\gamma(t)) - t).
$$
The equivalence between the notions of horosphere and level set of some Busemann function follows from Fact 7.

In Section 2 of the paper \cite{Busemann1932b}, Busemann establishes the following properties of horospheres:

\begin{enumerate}[i.)]
  \item Two horospheres associated to the same ray are equidistant.
  \item Horospheres are convex in the sense that every point not on a horosphere has a unique nearest 
  point on that horosphere.
  \item Horospheres in an $n$-dimensional line space are locally homeomorphic to standard spheres; therefore their    topological dimension is $n-1$.
\end{enumerate}

In Section 3, Busemann considers horospheres in a line space satisfying the limit circle axiom.
He first observes that the horospheres in such a space are flat subspaces. In particular they are 
themselves line spaces (of dimension one less).
Furthermore, any horosphere separates the  ambient space into two convex regions.

Busemann then specializes his investigation to dimension 3. He calls \emph{plane} a 
horosphere in  a 3-dimensional line space   satisfying the limit circle axiom. Any ray defining the horosphere is said to represent the \emph{normal direction}
of that plane.

Busemann then shows that there is a unique plane through three non aligned points 
(Facts 17 and 25).  
He also proves that given a point $P$ not on a plane there is a unique
footpoint (nearest point) for $P$ on the plane.

Note that this property implies that at any point on a sphere, there is a unique
plane tangent to that sphere (where by ``tangent'' one means that the
intersection of the plane and the sphere is a singelton).

In Section 4, Busemann discusses the geometry of 2-dimensional
line spaces satisfying the limit circle axiom.
He first observes that the horospheres in such a space are exactly the lines.
In the rest of the section he shows that such a space is homeomorphic to
$\mathbb{R}^2$ and Pash's axiom is satisfied. The arguments are similar to
those in \cite{Busemann1933}, see Section \ref{s:Pasch}.
It is also easy to show that in such a  space the parallel postulate (as formulated
in Hilbert's Axiom IV) holds.

In  Section 5, Busemann recapitulates the previous discussion in the following
claim: \emph{For a 3-dimensional line space 
satisfying the limit circle axiom, all the projective axioms of a 3-space have been verified}.
This means that all Hilbert's axioms of incidence and order, that is,  groups (I) and (II) are satisfied.
This also implies that Desargues' theorem holds in a suitable sense. Since the parallel postulate also holds,
$X$ can be mapped homeomorphically to $\mathbb{R}^3$ in such a way that the geodesics in 
this plane are mapped to straight lines in $\mathbb{R}^3$, see \cite{Troyanov} for a discussion of 
Desargues' theorem and its relation with such a mapping.

Busemann then observes that the spheres in $X$ are strictly convex and at each point they
admit a unique tangent plane. He then shows (Fact 39) that any pair of spheres in $X$ are related by a
similarity transformation, and from this he concludes that the line space is a Minkowski space (Fact 39).
This  completes  the proof of the main Theorem.

Busemann then proceeds to prove the converse, namely that  Minkowski metrics on $\mathbf{R}^3$
with strictly convex spheres and with a unique tangent plane at each point are line spaces satisfying
the limit circle axiom (Fact 41). He concludes  the paper by recalling that it is then sufficient to add a single axiom to conclude that the
line  space   is Euclidean. He states the additional axiom as a condition involving a configuration of
6 points proposed by O. Veblen. Note that many other equivalent conditions were proposed by a number of
mathematicians. The most elegant and famous among them is probably the \emph{parallelogram law} 
published in 1935 by P. Jordan and J. von Neumann in \cite{JordanvonNeumann}.

%_______________________
\section{Spaces with Convex Spheres and the Parallel Postulate}   \label{s:1933}

In the fifth postulate of Euclid's  \emph{Elements} and in Hilbert's axiom (IV), the notion of parallelism between two lines refers to the fact that the two lines are coplanar and disjoint. In the paper \emph{On spaces with convex spheres and the parallel postulate}
\cite{Busemann1933}, Busemann introduces a different notion of parallelism, which he will also adopt (with a small variation) in his later work, in particular, in his book \emph{The Geometry of geodesics}. 
The definition is the following:

\begin{definition}[Space satisfying the parallel postulate]   
A line space is said to satisfy the parallel postulate if whenever $g$ is a straight line, $P$ a point not on $g$, and 
$X_n$ a sequence of points on $g$ converging to infinity, the straight lines $PX_n$ always converge to the same 
straight line (depending only on $g$ and $P$.)
\end{definition}

Using this notion, Busemann gives a characterization of Minkowski spaces that uses spheres, rather than  horospheres as in the previous paper. More precisely, he first obtains the following result, valid in dimension $\geq 3$:

\begin{theorem} \label{th.Minkowski-31}
  If in a line space of dimension $\geq 3$ the distance spheres are convex and differentiable and the parallel axiom holds, 
then the space can be mapped topologically onto a Euclidean space in such a way that the geodesics are mapped to straight lines and the metric is Minkowskian.
\end{theorem}

The differentiability of a sphere, in Busemann's metric setting, means the existence and uniqueness of  ``as few tangent lines as possible''
at each point of the sphere.  This is made precise in \S 3 of his paper.

The case of dimension 2 requires different hypotheses. In Busemann's terminology, in a line space, the distance function is said to be convex if the distance function from an arbitrary fixed point to the points moving on an arbitrary line is convex. In the same manner, by taking the $\alpha$-th power of this function, one defines the notion of a line space in which the $\alpha$-th power of the distance function is convex.

In the case of dimension two, Busemann obtains the following.
\begin{theorem}\label{th.Minkowski-32}
In a two-dimensional line space in which the parallel postulate holds, if for some $\alpha\geq 1$ the $\alpha$-th power of the distance function is convex and differentiable then the same conclusion as in Theorem \ref{th.Minkowski-31} holds.
\end{theorem}

The outline of the paper is the following.

The paper  is divided into 5 sections.

 Section 1 concerns the notion of parallelism in line spaces. Busemann introduces the notion of asymptotic lines, and he makes a thorough investigation of this notion, making the relations with limit spheres and their orthogonal trajectories.  Busemann refers to the book \cite{Cartan} by \'E. Cartan in which the latter studied such a relation in Riemannian spaces of nonpositive curvature. In fact, the notion of asymptoticity of lines, in the setting of surfaces embedded in 3-space, was already studied by Hadamard in his paper \cite{Hadamard}. The study made by Busemann in his paper is a generalization of those made by Hadamard and Cartan to this setting of abstract line spaces. Using this notion of parallelism, Busemann states his \emph{parallel postulate} in the setting of line spaces. The definition involves 5  statements concerning asymptotic lines. Dimension 2 requires a special study, and for this case Busemann uses the result of his paper \cite{Busemann1933a}.

Section 2 concerns the convexity of spheres and its consequences. Busemann makes the relation of this notion with the existence of midpoints, and he notes that this was already noticed by Minkowski in his \emph{Geometrie der Zahlen} \cite{Minkowski-Zahlen}.  Busemann refers again to Cartan's book  \cite{Cartan} for the case of Riemannian manifolds. In this general setting,
 Busemann gives a formulation using the convexity of the distance function and the uniqueness of foots of perpendiculars. A thorough investigation of the notion of perpendicularity in this non-Riemannian setting is needed. In particular, the notion of perpendicularity is extended to a pair of lines which do not meet: a line $g$ is said to be perpendicular to a line $h$ if there is a parallel to $g$ which is perpendicular to $h$. Again, in the case of dimension two, a special study is made.  All these ideas will play an essential role in the later work of Busemann.

Section 3 concerns the characterization of Minkowski spaces of dimension $>2$. In this section, Busemann defines the notion of differentiablity of spheres in this metric setting. This involves the introduction of metric notions of  hypersurfaces, midsurfaces and hyperspheres, the use of the parallel postulate, and of several topological lemmas. 

In Section 4, Busemann studies the question of embedding a 2-dimensional line space into a 3-dimensional one. He recalls that Hilbert, in his \emph{Grundlagen}, already considered the question of the embedding of a plane satisfying all the axioms of order and intersections in a three-dimensional projective space. He showed that the embeddability is equivalent to the fact that Desargues property is satisfied. This is one of the starting points for Busemann's work on what he called Desarguesian spaces. In the paper under review, Busemann  adresses a similar question. The space is equipped with a metric and with a certain system of curves, called \emph{lines}.  Busemann proves that if in this space Desargues theorem holds, then the two-dimensional space can be embedded in a three-dimensional space such that these lines are sent to lines of the projective space. The proof uses a construction by Hessenberg \cite{Hessenberg}.

Section 5 concerns the characterization of 2-dimensional Minkowski spaces. The embedding of \S 4 is used. The convexity property of the $\alpha$-th power of the distance function to which we referred in Theorem \ref{th.Minkowski-32} above is equivalent to the fact that in that space, the spheres are convex. Busemann uses this fact to show the main theorem in dimension two, namely, if in a 2-dimensional line space the parallel postulate holds and if there exists $\alpha>1$ such that the $\alpha$-th power of all distance functions are convex and differentiable, then the metric is Minkowskian.

The paper ends with a discussion of the convexity of ellipses and hyperbolas, with their natural definition in a general 2-dimensional line space. This is naturally motivated by the result on the convexity of spheres. Busemann shows in particular that if in such a space the branches of any hyperbola are convex, then for any two points, the perpendicular bisectors exist and are lines. Then he proves that if in any two-dimensional line space the branches of any hyperbola are convex curves, then the metric is either Euclidean of hyperbolic.

Finally, let us note that some of the results of the three papers \cite{Busemann1932b,Busemann1933a, Busemann1933} were also reproved in the book \cite{Busemann1955}

\section*{Appendix: Hilbert's axioms}

In this section, we reproduce, for the readers' convenience, Hilbert's axioms from his \emph{Foundation of Geometry}. 
We are following the second English edition \cite{Hilbert-Grundlagen}, which is based on the tenth German edition. 
This edition differs marginally from the seventh edition, to which  Busemann refers.
Besides the axioms themselves,  Hilbert's treatise  also contains the consequences of the axioms as well as numerous comments enriching the axiomatic system. 

The goal of this system of axioms is to lay down foundations for Euclidean
geometry in 3-space. Three undefined notions (objects) are assumed: \emph{points} (denoted by capital Latin letters), \emph{lines} (denoted by small Latin letters) and \emph{planes}(denoted by small Greek letters). 
Some undefined relations are also assumed; these are the \emph{incidence} relations (a point is on a line, a line is contained in a plane etc.),
the \emph{order} relation  (a point lies between two others) and the \emph{congruence} relation formalizing the idea that two segments or
two angles can be superposed.

The axioms are divided into five groups describing the basic properties of incidence, order and congruence. The parallel postulate and 
a continuity (or completeness) axiom are also assumed.

\newpage

\textbf{I. Axioms of incidence}

\begin{enumerate}
\item[(I.1)]  For every two points $A, B$ there exists a line  that contains
each of the points $A,B$.

\item[(I.2)]  For any two distinct points $A, B$ there exists no more than one line  that contains
each of the points $A,B$.

\item[(I.3)]   There exist at least two points on a line. There exist at least three points that do not lie on a line.

\item[(I.4)]  For any three points  $A$, $B$, $C$  that do not lie on the same line there exits a plane $\alpha$ that
contains each of the points $A, B, C$. For every plane there exists a point which it contains.

\item[(I.5)]   For any three points $A$, $B$, $C$   that do not lie on one and the same line there exists no more than one plane that contains each of the three points $A, B, C$.

\item[(I.6)]  {If two points $A$, $B$ of a  line $a$ lie in a plane $\alpha$, then every point of $a$ lies in $\alpha$}.

\item[(I.7)]  {If two planes $\alpha$, $\beta$ have a point $A$ in common, then they have at least one more point $B$ in common.}

\item[(I.8)]  {There exist at least four points not lying in a plane.}
\end{enumerate}

\textbf{II. Axioms of order }

\begin{enumerate}
\item[(II.1)] 
If a point $B$ lies between a point $A$ and a point $C$ then the
points $A$, $B$, $C$  are three distinct points of a line, and  $B$ lies also between
$C$ and $A$.

\item[(II.2)]   For two points $A$ and $C$, there always exists at least one point
$B$ on the line $AC$ such that $C$ lies between $A$ and $B$.

\item[(II.3)]  {Of any three points on a  line, there exists  no more than one that lies between
the other two.}

\item[(II.4)]   Pasch's Axiom:  {Let $A$, $B$, $C$ be three points not lying in the
same  line and
let $a$ be a  line lying in the plane $ABC$ and not passing through any of the
points $A$, $B$, $C$. If  the  line $a$ passes through a point
of the segment $AB$, it   also passes through either a point of the segment
$BC$ or a point of the segment $AC$.}
\end{enumerate}

\medskip

\textbf{III. Axioms of congruence } 

\begin{enumerate}
\item[(III.1)]   {If $A$, $B$ are two points on a  line $a$,
and if $A'$ is a point on the same or another
 line $a'$, then it is  possible to find a point 
$B'$  on a given side of the line $a'$ so that the segment $AB$ 
is congruent to the segment $A'B'$. We  write  $AB \equiv A'B'$. }

\item[(III.2)]  If two segments are congruent to a third one they are congruent to each other.

\item[(III.3)] 
On the line $a$ let $AB$ and $BC$ be two segments which except for $B$ have no point in common. 
Furthermore, on  the same or  another line $a'$  let $A'B'$ and $B'C'$
be two segments having  no point other than $B'$ in
common.  In that case  if $AB \equiv A'B'$ and $BC \equiv B'C'$,
we have $AC \equiv A'C'$.
\end{enumerate}

The next two axioms   concern the congruence of planar angles
and are of little interest in Minkowski geometry.

\begin{enumerate}
\item[(III.4)]\footnote{We have summarized Axiom III.4 -- which is long in Hilbert's version -- while keeping the meaning.} 
  {Given an angle in a plane $\alpha$ and a half-line in another plane $\beta$ there exists exactly two half-lines in $\beta$
generating angles congruent with the given angle.}

\item[(III.5)]  {If $ABC$ and $A'B'C'$ are two triangles such that $A'B' \equiv AB$,  $A'C' \equiv AC$ and $\widehat{B'A'C'} \equiv \widehat{BAC}$,
then we also have $\widehat{ABC} \equiv  \widehat{A'B'C'}$ and $\widehat{ACB}\equiv \widehat{A'C'B'}$.}
\end{enumerate}

\newpage

\textbf{IV. The parallel postulate  } 

\begin{enumerate}
\item[(IV.)]
Let $a$ be any line and $A$ a point not on it. Then there is at most one line in the plane, determined by $a$ and $A$, that passes through $A$ and does not intersect $a$.
\end{enumerate}

\textbf{V. Axioms of continuity }

\begin{enumerate}
\item[(V.1)]  Archimedes' Axiom:
If $AB$ and $CD$ are any segments then there exists a number $n$ such that $n$
 segments $CD$ constructed contiguously from $A$, along the ray from $A$ through $B$, will pass beyond the point $B$.
 
\item[(V.2)]  Completeness' Axiom: 
An extension of a set of points on a line with its order and congruence relations that would preserve 
the relations existing among the original elements as well as the fundamental properties of
 line order and congruence that follows from Axioms I-III, and from V.1 is impossible.
 \end{enumerate}

\subsection*{On the triangle inequality} 

Before ending this paper, we aim at clarifying a point. In his comment to his IVth problem, Hilbert suggests to consider the
triangle inequality as an \emph{Axiom} and not as a \emph{Theorem} (see Hilbert's quotation on page \pageref{quote.hilbert} 
in the present paper). However the triangle inequality is \emph{not} explicitly stated as a theorem in Hilbert's foundations.

To understand in what sense Hilbert considers the triangle inequality as a theorem, let us recall that in Euclid's Elements, 
Proposition 20 of Book 1 says that \emph{In any triangle the sum of any two sides is greater than the remaining one.}
The proof derives from Proposition 5: \emph{On an isosceles triangle the angles at the base are equal to one another}  and
Proposition 19:  \emph{In any triangle the greater angle is subtended by the greater side}.

In his \emph{Foundations}, more precisely, in the comments following the congruence axioms, Hilbert  discusses the notions of angle and segment comparison. Axiom (III.3) allows also an obvious definition of 
segment addition. With this material, Hilbert proves the equality of the base angles in an isosceles triangle as Theorem 11 and he proves that 
in every triangle the greater angle lies opposite to the greater side as 
Theorem 23. One may then  prove the triangle inequality exactly as in Euclid.
This justifies Hilbert's claim that it is a  \emph{proposition} that \emph{ in every triangle the sum of two sides is greater than the third}
and gives a meaning to his suggestion to assume this fact as an \emph{axiom} instead.
 
Let us end with two short comments. First we stress that in the above proof the angle congruence Axiom plays an important role. Thus the triangle inequality in classical Euclidean geometry appears as a consequence of some properties of angles. Finally let us mention that although Hilbert did not explicitly state the triangle inequality in the \emph{Foundations}, he wrote a paper  in 1902 about the role of the theorem on the base angles in an isosceles triangle \cite{Hilbert1902}. This paper is reproduced in appendix II of  \cite{Hilbert-Grundlagen}. There, the triangle inequality is explicitly stated as a proposition.

\medskip

\noindent {\bf Acknowledgement} The authors would like to thank Valerii Berestovskii for reading this paper and sending corrections. The first author acknowledges support from the ANR program FINSLER.

%_______________________________ 


\begin{thebibliography}{999}

\bibitem{Alexandroff} P.~Alexandroff, Dimensionstheorie. Ein Beitrag zur Geometrie der abgeschlossenen Mengen, Math. Annalen 106 (1932), p.~161-238.

\bibitem{Aristotle-Metaphysics} Aristotle, the Metaphysics, In \emph{The Complete Works of Aristotle: The Revised Oxford Translation} (J.  Barnes, editor), Volume 2, 1552-1728, Translated by W. D. Ross, Princeton University Press, Princeton, 1984. 


\bibitem{Berestovskii} V. N. Berestovskii, Busemann's results, ideas, questions and locally compact
homogeneous geodesic spaces, this volume.



\bibitem{Bonola} R. Bonola,
\textit{La geometria non-euclidea, Esposizione storico-critica del suo sviluppo}.  First edition, Ditta Nicola Zanchinelli editore, 
Bologna 1906. German translation by M. Liebmann in the collection {\it Wissenschaft und Hypothese}, Teubner, Leipzig, 1908. English
  translation by H. S. Carslaw: \textit{Non-Euclidean Geometry, A critical and historical study of its development}. First edition, 1912.  Reprinted by Dover, 1955.
  

\bibitem{Busemann1932b}  H. Busemann,
\emph{\"Uber die Geometrien, in denen die  ``Kreise mit unendlichem Radius''   die k\"urzesten Linien sind.} 
Math. Ann. 106 (1932), no. 1, 140--160. English translation by A. A'Campo:  On the geometries where circles of infinite radius are the shortest lines, this volume.


 \bibitem{Busemann1933a} H. Busemann,   \emph{Paschsches Axiom und Zweidimensionalit\"at}. 
Math. Ann. 107 (1933), no. 1, 324--328.  English translation by A. A'Campo:  
Pasch's Axiom and Two--Dimensionality, this volume.



\bibitem{Busemann1933} H. Busemann, 
\emph{\"Uber  R\"aume mit konvexen Kugeln und Parallelenaxiom. }
Nachrichten G\"ottingen 1933, 116-140 (1933).  English translation by A. A'Campo:  
 On spaces with convex spheres and the parallel postulate, this volume.
 
 

\bibitem{Busemann1942} H. Busemann,   
\emph{Metric Methods in Finsler Spaces and in the Foundations of Geometry.}
Annals of Mathematics Studies, no. 8. Princeton University Press, Princeton, N. J., 1942


 \bibitem{Busemann-Minkowski1950}  H. Busemann, The foundations of Minkowskian geometry. Comment. Math. Helv. 24, (1950). 156-187

\bibitem{Busemann1955} H. Busemann,   {\it The geometry of geodesics}, Academic Press  (1955), reprinted by Dover in 2005.


\bibitem{B-Minkowski1969}  H. Busemann,  Transitive geodesics in Minkowski planes. Math. Chronicle 1 1969/1970 part 1, 27-29.


\bibitem{Busemann1970} H. Busemann,  
Recent synthetic differential geometry. Berlin-Heidelberg-New York: Springer-Verlag 1970.


\bibitem{BP1979} H. Busemann and B. B. Phadke,  Minkowskian geometry, convexity conditions and the parallel axiom. J. Geom. 12 (1979), no. 1, 17-33.


%
% \bibitem{BW}  
%Bryant, V. W.; Webster, R. J.
%Convexity spaces. I. The basic properties. 
%J. Math. Anal. Appl. 37 (1972),  206--213.  

%
% \bibitem{Doignon} Doignon, J.-P.  "Caract\'erisations d'espaces de Pasch-Peano".
%Acad. Roy. Belg. Bull. Cl. Sci. (5) 62 (1976), no. 9, 679--699.
%
% \bibitem{Cantwell} Cantwell, John; Kay, David C. Geometric convexity. III. Embedding. Trans. Amer. Math. Soc. 246 (1978),  211--213

\bibitem{Cartan} \'E. Cartan, Le\c cons sur la g\'eom\'etrie des espaces de Riemann, Paris, Gauthier-Villars, 1928.

% \bibitem{Coppel} W. A. Coppel, Foundations of Convex Geometry. Australian Mathematical Society Lecture Series, 12. (Cambridge University Press, 1998).  

% \bibitem{Coxeter} H. S. M. Coxeter, Introduction to Geometry. Wiley Classics Library (Wiley, 1989).
 
%\bibitem{Pasch}  M. Pasch, Vorlesungen  ueber neuere Geometrie. 2. Auflage (Springer-Verlag, 1926).


 \bibitem{Frechet} M. Fr\'echet,   Sur quelques points de calcul fonctionnel, th\`ese de doctorat publ., dans Rendiconti del Circolo matematico di Palermo, t. 22 (1906), 1--74.
 
\bibitem{Goenner} H. Goenner, Documents concerning Busemann, from the library of G\"ottingen University, this volume.
 

 
 \bibitem{Hadamard} J. Hadamard, Les surfaces \`a courbures oppos\'ees et leurs lignes g\'eod\'esiques. Journal de math\'ematiques pures et appliqu\'ees 5e s\'erie, tome 4 (1898), p. 27--74.
 
 \bibitem{Hessenberg} G. Hessenberg, Grundlagen der Geometrie, Herausgegeben von W. Schwan, Berlin and Leipzig, Verlag Walter de Gruyter \& Co. 1930.
 
 
\bibitem{Hilbert1902} D. Hilbert,
Ueber den Satz von der Gleichheit der Basiswinkel.
Proceedings of the London Math. Soc., Vol. 35. 1902
 
\bibitem{Hilbert-Problems}  D. Hilbert,  Mathematische Probleme, \emph{G\"ottinger Nachrichten}, 1900, pp. 253--297, reprinted in {\it Archiv der Mathematik und Physik}, 3d. ser., vol. 1 (1901) pp. 44--63 and 213--237. English version, ``Mathematical problems", translated by M. Winston Newson, \emph{Bulletin of the AMS}, vol. 8, 1902, pp. 437-- 445 and 478--479. Reprinted in
 the  Bull. Amer. Math. Soc. (N.S.) 37 (2000), no. 4, 407-436. French edition, \emph{Sur les probl\`emes futurs des math\'ematiques}, 1902, trad. L. Laugel.


\bibitem{Hilbert-Grundlagen} D. Hilbert, 
Foundations of geometry. 
Second edition. Translated from the tenth German edition by Leo Unger Open Court, LaSalle, Ill. 1971. 

%
%\bibitem{Hilbert2004}
%\emph{David Hilbert's lectures on the foundations of geometry, 1891-1902. }
%Edited by Michael Hallett and Ulrich Majer. David Hilbert's Foundational Lectures, 1. 
%Springer-Verlag, Berlin, 2004.

 \bibitem{JordanvonNeumann}  
P. Jordan and J. von Neumann, 
On Inner Product in Linear Metric Spaces, 
Ann. of Math. 36 (1935), no. 2, 719-723.


 \bibitem{Loba-Geometrische} N. I. Lobachevsky,
 \emph{Geometrische Untersuchungen zur theorie der Parallellinien} (Geometrical researches on the theory of parallels). Kaiserl. russ. wirkl. Staatsrathe und ord. Prof. der Mathematik bei der Universit\"at Kasan.  Berlin 1840. In der Finckeschen Buchhandlung. 61 pages.
 French translation by J. Ho\"uel,  Etudes g\'eom\'etriques sur la th\'eorie des parall\`eles, \textit{M\'emoires de la Soci\'et\'e de sciences physiques et naturelles de Bordeaux}. t. IV, 1866, pp. 83--182, and Paris, Gauthier-Villars, 1866. 
 English translation by G. B. Halsted, Geometrical researches on the theory of parallels, first published by the University of Texas at Austin, 1891. Reprinted in Bonola \cite{Bonola}.   
 
 \bibitem{Pangeometry} N. I. Lobachevsky, Pang\'eom\'etrie ou pr\'ecis de g\'eom\'etrie fond\'ee sur une th\'eorie g\'en\'erale
et rigoureuse des parall\`eles.  
In: Collection of memoirs written by professors of the University of Kazan on the occasion
of the 50th anniversary of its foundation, Vol. I, Kazan, 1856, 277--340. English edition: Pangeometry,  Translation, notes and commentary by A. Papadopoulos, Heritage of European Mathematics, Vol. 4, European Mathematics Publishing House,  322 pages, 2010.



 
 \bibitem{Menger1928} K. Menger, 
Untersuchungen \"uber allgemeine Metrik.  
Math. Ann. 100, 75-163 (1928).
 
 
 
% 
% \bibitem{Menger1930} K. Menger, 
%Untersuchungen \"uber allgemeine Metrik IV.
%Math. Ann. 103 (1930), no. 1, 466--501. 
%
%
% \bibitem{Menger1931} K. Menger,  
%New Foundation of Euclidean Geometry. 
%Amer. J. Math. 53 (1931), no. 4, 721--745.  

\bibitem{Minkowski-Zahlen} H. Minkowski, Geometrie dre Zahlen, Leipzig, 1896. (Busemann refers to the 1910 edition.)

\bibitem{Papadopoulos2005} 
A. Papadopoulos,  
Metric spaces, convexity and nonpositive curvature. 
IRMA Lectures in Mathematics and Theoretical Physics, 6. 
European Mathematical Society (EMS),  Zurich, 2005, Second edition, 2014.

% \bibitem{Peano} G. Peano, I principii di geometria logicamente esposti (Fratelli Bocca editori, 1889) (also G. Peano, Opere scelte, vol. II, Edizione Cremonese, Roma, 1958, pp. 56--91.    


  \bibitem{2012-Hilbert1} A. Papadopoulos, Hilbert's fourth problem. In: Handbook of Hilbert geometry (A. Papadopoulos and M. Troyanov, ed.), European Mathematical Society Publishing House,  IRMA Lectures in Mathematics and Theoretical Physics, Vol. 22, p. 391-431, 2014.

\bibitem{RH} R. Rashed and Ch. Houzel, Th\=abit et la th\'eorie des parall\`eles, Arabic Sciences and Philosophy, 15.1, 2005, p. 9--55.


\bibitem{R}  R. Rashed, Le mouvement en g\'eom\'etrie classique,  Al-Mukhatabat (2013), p. 58-68.


\bibitem{Troyanov} M. Troyanov,  On Pasch's Axiom and Desargues' Theorem  in Busemann's work, this volume.
 
% \bibitem{Veblen} O. Veblen, A system of axiom for geometry. Trans. Amer. Math. Soc. 5, 432--384
%  (1904). 
 
\end{thebibliography}
\end{document}